\pdfoutput=1 
\documentclass[12pt,draftcls,onecolumn]{IEEEtran}

\usepackage{graphicx}
\usepackage{subfigure}
\usepackage{amsmath}
\usepackage{amsfonts}
\usepackage[psamsfonts]{amssymb}
\usepackage{comment}
\usepackage{times}
\usepackage{changepage}
\usepackage{cases}
\usepackage{mathtools}
\usepackage{color}
\usepackage{booktabs}
\usepackage{array}
\usepackage{mathrsfs}
\usepackage{enumerate}
\usepackage{setspace}
\usepackage[mathscr]{eucal}
\usepackage[scr=euler,calscaled=.94,bb=ams,frak=mma,frakscaled=.97,scr=rsfs]{mathalfa}
\mathtoolsset{showonlyrefs=true}

\newtheorem{definition}{Definition}

\newtheorem{lemma}{Lemma}

\newtheorem{theorem}{Theorem}
\newtheorem{algorithm}{Algorithm}
\newtheorem{assumption}{Assumption}

\newtheorem{problem}{Problem}

\newenvironment{customprob}[1]
  {\renewcommand\theproblem{#1}\problem}
{\endproblem}

\newcommand{\argmin}{\mathop{\rm arg~min}\limits}

\newcommand{\0}{\emptyset}
\newcommand{\R}{\mathbb{R}}
\newcommand{\N}{\mathbb{N}}

\newcommand{\bmu}{\boldsymbol{\mu}}

\newcommand{\Lap}{\textnormal{Lap}}

\newcommand{\tl}{\tilde{\ell}}
\newcommand{\bx}{\boldsymbol{x}}
\newcommand{\lqr}{\tl_{q}^{r}}
\newcommand{\lsp}{\tl_p^s}

\newcommand{\Adj}{\textnormal{Adj}}
\renewcommand{\P}{\mathbb{P}}
\newcommand{\E}{\mathbb{E}}

\newcommand{\new}[1]{#1}
\newcommand{\blue}[1]{\new{#1}}

\let\originalleft\left
\let\originalright\right
\renewcommand{\left}{\mathopen{}\mathclose\bgroup\originalleft}
\renewcommand{\right}{\aftergroup\egroup\originalright}

\title{\LARGE \bf
\blue{Cloud-Enabled Differentially Private Multi-Agent Optimization with Constraints}
}

\author{Matthew T. Hale and Magnus Egerstedt\thanks{The authors are with the School of Electrical and Computer Engineering, Georgia Institute of
Technology, Atlanta, GA 30332, USA. Email: \texttt{\{matthale, magnus\}@gatech.edu}. Research supported in
part by  the NSF under Grant CNS-1239225.}
}

\begin{document}

\maketitle
\thispagestyle{empty}
\pagestyle{empty}

\begin{abstract}
We present an optimization framework for solving
multi-agent nonlinear programs subject to inequality constraints while
keeping the agents' \new{state trajectories} 
private. Each agent has an objective
function depending only upon its own state and the agents
are collectively subject to global constraints. 
The agents do not directly communicate with each other but instead
route messages through a trusted cloud computer. 
The cloud computer
adds noise to data being sent to the agents in accordance with the framework of
differential privacy in order to keep each agent's \new{state trajectory} private from all other agents and
any eavesdroppers. This private problem can be viewed
as a stochastic variational inequality and is solved 
using a projection-based method for solving variational inequalities
that resembles a noisy primal-dual gradient algorithm.
Convergence of the optimization algorithm
in the presence of noise is proven and a quantifiable
trade-off between privacy and convergence is extracted
from this proof. 
Simulation results
are provided that demonstrate numerical convergence for both 
$\epsilon$-differential privacy and $(\epsilon,\delta)$-differential
privacy.
\end{abstract}

\section{Introduction}
Optimization problems spread across teams of agents arise naturally in several
fields, including communications \cite{kelly98,chiang07},
robotics \cite{cortes02,soltero13}, machine learning \cite{boyd11},
sensor networks \cite{khan09,trigoni12}, and 
smart power grids \cite{vytelingum10,caron10,nazari14}. Correspondingly,
a variety of approaches have been developed that solve problems
with a wide variety 
of formulations. For example, \cite{nedich09} allows for
distributed optimization of non-differentiable objectives with
time-varying communication links, \cite{lobel11} considers a similar
problem formulation in which communication links fail over time, and
\cite{wei13} uses a distributed Newton method to solve dynamic
network utility maximization problems. Many other problem types
and solution schemes exist in the literature, and a broad 
exposition of results can be found in \cite{bertsekas97}. 

\new{
In some cases multi-agent optimization is done using sensitive user data.
A concrete example of such a case comes from smart power grids. In smart power grids,
homeowners share their power usage information with others on the grid
to allow network management (e.g., frequency regulation \cite{nazari14}) 
and to minimize their own power costs. In some cases, the granular
power usage data shared in smart grids can be used to infer sensitive details
of users' personal lives \cite{doe10,edps12}. In particular,
smart grid data can ``provide a detailed breakdown
of energy usage over a long period of time, which can show patterns of use,''
\cite[Page 15, Item 16]{edps12}. Further, given these patterns, ``[p]rofiles
can thus be developed and then applied back to individual households
and individual members of these households,'' \cite[Page 15, Item 18]{edps12}.
These usage patterns in turn ''could reveal personal details about the lives
of consumers, such as their daily schedules,'' \cite[Page 2, Paragraph 5]{doe10}.

It is precisely the deduction of
such patterns that we wish to prevent in the context of multi-agent optimization.
Based on the potentially revealing nature of some user data,
we seek to optimize while protecting sensitive user data both from
eavesdroppers and other agents in the network.
In some sense, privacy and optimization are competing objectives in that
agents who only seek to optimize may freely share their states with others
in the network, while agents concerned only with privacy may be inclined
to share no information at all. To privately optimize, then, we must strike
a balance between these two different, competing objectives. 
}

One approach to privacy that has recently seen widespread use is
\emph{differential privacy}. Differential privacy was originally 
established in the database literature and keeps
\new{sensitive} database entries private when a database is queried by adding noise
to the result of that query \cite{dwork06a,dwork06c}.
The authors of \cite{dwork08,dwork13} survey some of the
important developments in this vein.
Differential privacy has been adapted to dynamical systems in order to
keep \new{sensitive} inputs private from an adversary observing a system's outputs \cite{pappas14}. 
A dynamical system is differentially private if
inputs that are close in the input space produce outputs that have similar probability distributions; 
these notions will be made precise in Section~\ref{sec:privacy}.

\new{It is the dynamical systems notion of differential privacy that
we apply to keep agents' state trajectories private while optimizing.} 
One appealing aspect of differential privacy is its resilience to
post-processing, which allows for arbitrary processing of 
private information without the threat of its privacy guarantees being weakened
\cite[Theorem 1]{pappas14}. Differential privacy is also robust
to arbitrary side information, meaning that 
an adversary cannot weaken
differential privacy by much through using information gleaned from another source \cite{kasiv08}.


There has already been some work on enforcing differential privacy in optimization. 
In \cite{hsu14} linear programs are solved in a framework that
allows for keeping objective functions or constraints private. 
The authors of \cite{han14} consider a similar setting
wherein linearly constrained problems with affine objectives
are solved while keeping the objective functions private.
In the multi-agent setting, \cite{huang15} solves distributed
consensus-type problems while keeping the agents' objective
functions private, while \cite{huang12} solves similar problems
while keeping each agent's initial state private. 

In this paper we solve non-linear programs wherein
\new{
each agent's state trajectory is sensitive information and the agents
therefore seek to protect their exact state trajectories from other
agents and any eavesdroppers.
}
\new{To protect these sensitive data}, a trusted cloud computer
is used that performs certain computations upon information
it receives from the agents, makes the results of those
computations private by adding noise to them, and then sends the 
private results to each agent. Each agent then updates its
state locally using the information it received from the cloud,
and this process of sharing and updating information is repeated. 

\new{
Our motivation for developing a mixed centralized/decentralized algorithm
is inspired by the prominence of cloud computing in many real-world applications.
A survey of existing cloud applications is given in \cite{zhang10}, and 
that reference elaborates on the scalability of the cloud and its ability
to coordinate many mobile devices. It is precisely
these features of the cloud that make it an attractive choice here. 
In this paper, the cloud,
viewed as a central aggregator, is an integral part of the optimization process, and
we leverage its scalability to aggregate ensemble-level information, perform computations
upon that data in a private manner, and then distribute these private
results to the agents. 
}

\new{
The privacy implementation in this paper differs from the aforementioned
references on private optimization in several key ways. We are interested in solving problems in
which the agents collectively run an on-line optimization algorithm collaboratively by
sharing (private functions of) sensitive information. In the problems we consider, each iteration of 
the optimization algorithm determines each agent's next state. That is, the
iterates of the optimization algorithm \emph{are} the agents' states, and
it is each agent's desire to keep its state trajectory private to protect
information about its behavior. 
Accordingly, while the above references on private
optimization keep other problem data private,
here we must keep entire trajectories of states private while optimizing. 
In addition, we incorporate both nonlinear inequality constraints and
set constraints, which, to our knowledge, has not been explored in
other privacy implementations. 
}

\new{
Given the need to optimize while remaining private,
encryption alone cannot provide the 
privacy guarantees that are needed in the problems we examine.
In the ``upstream'' direction, encryption could be used to protect communications sent from
the agents to the cloud, provided the cloud could decrypt them.
However in the ``downstream'' direction, when the cloud sends transmissions to the agents, any encrypted
messages from the cloud would naturally need to be decrypted by the agents to allow each agent
to update its state. 
While this strategy can protect transmissions of sensitive data from eavesdroppers, 
having the agents decrypt transmissions from the cloud would 
expose all agents' sensitive data to each agent in the network, violating
the privacy guarantees that are required by each agent. 
Instead, what is
required here is a privacy implementation that protects user data from
eavesdroppers \emph{and} all others in the network, while still making that data useful
for optimizing. 
It is for this reason that we use differential privacy.
}

A preliminary version of this work appeared in \cite{hale15}.
The current paper adds a proof of convergence, a convergence estimate, 
\new{quantifies
the privacy-convergence trade-off}, and provides new
numerical results for two different privacy mechanisms.  
The rest of the paper is organized as follows. Section
II lays out the problem to be solved and its method
of solution. Next, Section III covers the necessary
elements of differential privacy and relates
them to the setting of optimization. Then
Section IV provides a proof of convergence for the
optimization algorithm used here and a bound
on its convergence, \new{in addition to exploring
the trade-off between privacy and convergence.}
Next, Section V
provides simulation results to support the theoretical
developments made. Finally, Section VI concludes
the paper.

\section{Optimization Problem Formulation} \label{sec:problem}
In this section we lay out the problem to be solved. First in 
Section~\ref{ss:problem} we lay
out the multi-agent problem and then, to aid in the exposition
of its solution method, formulate an equivalent
ensemble problem. Then the solution to that problem will be discussed
and, in Section~\ref{ss:comms},
will be adapted to the cloud-based architecture used here. 

\vspace{-0.1cm}

\subsection{Problem Overview} \label{ss:problem}
Consider $N$ agents indexed over the set $I := \{1, \ldots, N\}$, with
agent $i$ having state $x_i \in \R^{n_i}$ for some $n_i \in \N$. 
Agent $i$ seeks to minimize the objective function
\begin{equation}
f_i : \R^{n_i} \to \R,
\end{equation}
where $f_i$ depends only upon $x_i$, that is, each agent's 
objective function has no dependence upon the other agents' states. 
Using the notation $\nabla_{i} f_i := \frac{\partial f_i}{\partial x_i}$,
we state the following assumption for objective functions.
\begin{assumption} \label{as:objectives}
The function $f_i$ is $C^1$ and convex, and $\nabla_{i} f_i$ is Lipschitz with constant $L_i$ for all $i \in I$. \hfill $\triangle$
\end{assumption}

Assumption~\ref{as:objectives} allows for a broad class of functions to be used as objective functions, including any $C^2$ convex function
on a compact, convex domain (cf. Assumption~\ref{as:set} below).
Each agent's state is constrained to lie in a given set which we express as 
\begin{equation}
x_i \in X_i \subset \R^{n_i}. 
\end{equation}
Regarding each set $X_i$, we state the following assumption.
\begin{assumption} \label{as:set}
Each set $X_i$ is non-empty, compact, and convex.  \hfill $\triangle$
\end{assumption}

In particular Assumption~\ref{as:set} admits box constraints which are common in some multi-agent problems. 

Now define the ensemble state vector 
\begin{equation}
x = \left(\begin{array}{c} x_1 \\ x_2 \\ \vdots \\ x_N \end{array}\right) \in \R^n,
\end{equation}
where $n = \sum_{i=1}^{N} n_i$. We impose global inequality constraints on the agents
by requiring
\begin{equation}
g(x) := \left(\begin{array}{c} g_1(x) \\ g_2(x) \\ \vdots \\ g_m(x) \end{array}\right) \leq 0,
\end{equation}
where the above inequality is enforced component-wise, i.e., $g_j(x) \leq 0$
for all $j \in J := \{1, \ldots, m\}$. 
We now state our assumptions on $g$. 
\begin{assumption} \label{as:constraints}
The function $g:\R^n\to \R^m$ is $C^1$ and convex. In addition, 
for $\R^m$ and $\R^n$ both equipped with the same $p$-norm, 
the function $g_{x_j} := \frac{\partial g}{\partial x_j}$ is Lipschitz continuous with constant $K^j_{p}$ for all $j \in J$ with
respect to the metric induced by the $p$-norm. In addition
the function $g$ is Lipschitz with constant $K^g_{p}$ with respect to the same metric. \hfill $\triangle$
\end{assumption} 

In this paper we focus on the cases of $p = 1$ and $p = 2$. 
Like Assumption~\ref{as:objectives}, Assumption~\ref{as:constraints} allows for any convex, $C^2$ functions to be used
for constraints whenever Assumption~\ref{as:set} holds. 
We also have the following assumption on $g$. 

\begin{assumption} \label{as:slater}
The constraints satisfy Slater's condition, namely there exists a point $\bar{x} \in X$ such that $g(\bar{x}) < 0$. \hfill $\triangle$
\end{assumption} 

Assumption~\ref{as:slater} is commonly enforced in nonlinear programming problems to guarantee that strong 
duality holds. 
\new{
Under Assumptions~\ref{as:objectives}-\ref{as:slater}, we state an ensemble-level optimization
problem. To do so, we define the ensemble objective
\begin{equation}
f(x) = \sum_{i=1}^{n} f_i(x_i),
\end{equation}
and the set
\begin{equation}
X = \prod_{i=1}^{n} X_i,
\end{equation}
where the product is meant in the Cartesian sense. 
To fix ideas, we state the following optimization
problem that does not yet incorporate privacy; privacy will be formally included
in Problem~\ref{prob:withprivacy} in Section~\ref{sec:privacy}. 

\begin{customprob}{0.1} \label{prob:global} \new{\emph{(Preliminary; no privacy requirement)}}
\begin{align*}
\textnormal{minimize } &f(x) \\
\textnormal{subject to } &g(x) \leq 0\\
&x \in X.
\end{align*}
\end{customprob} \hfill $\lozenge$
}

\new{
We note here that Problem~\ref{prob:global} will be solved without having agent $i$ share
$f_i$ or $X_i$ with the other agents or with the cloud because these data are considered
sensitive information. Similarly, $g$ is considered sensitive and the cloud does
not share $g$ with any of the agents. 
}
The Lagrangian associated with Problem~\ref{prob:global} is
\begin{equation}
L(x, \mu) = f(x) + \mu^Tg(x),
\end{equation}
where $\mu$ is a vector of Kuhn-Tucker multipliers in the non-negative orthant of 
$\R^m$, denoted $\R^m_{+}$. 
Under Assumptions~\ref{as:objectives}, \ref{as:set}, and \ref{as:constraints}
a primal solution $\hat{x}$ exists and the set of all primal solutions is non-empty and compact. 
With the addition of Assumption~\ref{as:slater}, a dual solution $\hat{\mu}$ exists and the optimal primal and dual
values are equal \cite[Proposition 6.4.3]{bertsekas03}. 

Under Assumptions 1-4, a point $\hat{x}$ solves Problem~\ref{prob:global} if
and only if there exists a point $\hat{\mu} \in \R^m_{+}$ such that $(\hat{x}, \hat{\mu})$ is a saddle point of
$L$, that is, if and only if the point $(\hat{x}, \hat{\mu})$ satisfies
\begin{equation} \label{eq:saddle}
L(\hat{x}, \mu) \leq L(\hat{x}, \hat{\mu}) \leq L(x, \hat{\mu})
\end{equation}
for all $(x, \mu) \in X \times \R^m_{+}$ \cite[Proposition 6.2.4]{bertsekas03}. It is as saddle points of 
$L$ that we seek solutions $(\hat{x}, \hat{\mu})$ to Problem~\ref{prob:global}. 



Toward that end, we next define the symbols
\begin{equation}
L_{x} := \frac{\partial L}{\partial x}, \,\, L_{\mu} := \frac{\partial L}{\partial \mu}, g_x := \frac{\partial g}{\partial x}, \textnormal{ and } f_x := \frac{\partial f}{\partial x},
\end{equation}
and define the map
\begin{equation}
G(x, \mu) = \left(\begin{array}{r} L_{x}(x, \mu) \\ -L_{\mu}(x, \mu) \end{array}\right).
\end{equation}
In what follows, it is necessary for $G$ to be a Lipschitz mapping. Though the maps $f_x$ and
$g_x$ are Lipschitz by Assumptions~\ref{as:objectives} and \ref{as:constraints}, $G$ itself
cannot be shown to be Lipschitz because its domain, $X \times \R^m_{+}$, is unbounded by virtue
of $\R^m_{+}$ being unbounded. To rectify this situation, we use
Equation \eqref{eq:saddle}
to find a non-empty, convex, compact set containing $\hat{\mu}$ as was done in \cite{uzawa58}. 
From the second inequality in Equation \eqref{eq:saddle} we have
\begin{equation}
f(\hat{x}) + \hat{\mu}^Tg(\hat{x}) \leq f(\bar{x}) + \hat{\mu}^Tg(\bar{x})
\end{equation}
for $\bar{x}$ the Slater point as defined in Assumption~\ref{as:slater}. By the complementary slackness condition
we have
\begin{equation}
f(\hat{x}) \leq f(\bar{x}) + \hat{\mu}^Tg(\bar{x}). 
\end{equation}
Rearranging we find
\begin{equation}
\sum_{j=1}^{m} \hat{\mu}_j \leq \frac{f(\bar{x}) - f(\hat{x})}{\min\limits_{1 \leq j \leq m} \left\{-g_j(\bar{x})\right\}} \leq \frac{f(\bar{x}) - f(x^*)}{\min\limits_{1 \leq j \leq m} \left\{-g_j(\bar{x})\right\}},
\end{equation}
where $x^* \in \argmin_{x \in X} f(x)$. 
We then define the set 
\begin{equation}
\mathbb{M} := \left\{\mu \in \R^m_{+} : \|\mu\|_{1} \leq \frac{f(\bar{x}) - f(x^*)}{\min\limits_{1 \leq j \leq m} \left\{-g_j(\bar{x})\right\}}\right\}
\end{equation}
which is non-empty, compact, and convex by definition, and which contains $\hat{\mu}$. 
For economy of notation, we define the symbols $Z := X \times \mathbb{M}$ and 
$\hat{z} := (\hat{x}, \hat{\mu})$, and we will use $z := (x, \mu)$ to denote an arbitrary
point in $Z$.


Since $L(\cdot, \mu)$ is convex for all $\mu \in \mathbb{M}$ and $L(x, \cdot)$ is concave for all $x \in X$, 
we see that $G$ is monotone  \cite[Theorem A]{rockafellar70}.
Under Assumptions 1-4, a primal-dual pair $(\hat{x}, \hat{\mu})$ is a saddle point of $L$
if and only if it solves the following variational inequality (VI) \cite[Corollary 11.1]{konnov07}.
\begin{customprob}{0.2} \label{prob:vi} \new{\emph{(VI formulation; no privacy requirement)}}
Find a point $\hat{z} \in Z$ such that
\begin{equation}
\langle z - \hat{z}, G(\hat{z}) \rangle \geq 0
\end{equation}
for all $z \in Z$.\hfill $\lozenge$
\end{customprob} 
\new{
Further discussion on the equivalence of Problems~\ref{prob:global} and \ref{prob:vi} is given
in \cite{facchinei03}, Sections 1.3.1, 1.3.2, and 1.4.1.
Privacy is formally added to Problem~\ref{prob:vi} in the statement of Problem~\ref{prob:withprivacy}
in Section~\ref{sec:privacy}. 
}

We will use the notation $VI(K, F)$ to denote the generic problem of finding a point $x \in K$ such that
\begin{equation}
\langle y - x, F(x) \rangle \geq 0
\end{equation}
for all $y \in K$, and we will use the notation $SOL(K, F)$ to denote the solution set of $VI(K, F)$. 
The symbols $Z$ and $G$ refer to the specific problem under consideration in this paper
 so that Problem~\ref{prob:vi} is denoted $VI(Z, G)$ and its solution set is $SOL(Z, G)$. 
It is in the setting of variational inequalities that we will proceed and we 
focus on solving Problem~\ref{prob:vi} with the understanding that its solutions also solve Problem~\ref{prob:global}.

For a compact set $K$ and a monotone map $F$, one method of solving the variational inequality $VI(K, F)$ is using a projection method with an iterative Tikhonov regularization
as was done for deterministic variational inequalities in \cite{bakushinskii74} and 
for stochastic variational inequalities in \cite{poljak78}; these methods regularize the earlier Goldstein-Levitin-Polyak method
for solving such problems \cite{goldstein64,levitin61}.
The basic principle underlying these methods is that a point in $SOL(K, F)$ can be approached iteratively
with $F$ specifying the direction in which to move at each iteration. To endow this procedure with greater
numerical stability and, as will be shown, robustness to noise, the $k^{th}$ iteration specified in
\cite{bakushinskii74,poljak78} instead uses the direction specified by $F + \alpha_kI$ with $I$ the identity map, $\alpha_k > 0$, and $\alpha_k \to 0$. 
When $F$ is monotone, each map $F + \alpha_kI$ is strongly monotone so that $SOL(K, F + \alpha_kI)$ is a singleton. 
Letting $\xi_k$ denote the (unique) element of $SOL(K, F + \alpha_kI)$, for $\alpha_k > 0$ and $\alpha_k \to 0$
we have $\xi_k \to y_0$ where $y_0$ is the least-norm element of $SOL(K, F)$ (which itself is non-empty because $K$ is compact and $F$ is monotone). 
%

Given an initial point $z(0) \in Z$, the deterministic form of the regularized method to solve $VI(Z, G)$ is 
given below in Algorithm~\ref{alg:tikh}.
\begin{algorithm} \label{alg:tikh}
Given a point $z(0) \in Z$, apply the update law
\begin{align}
z(&k+1) = \Pi_{Z}\left[z(k) - \gamma_k\left(G\big(z(k)\big) + \alpha_kz(k)\right)\right] \\
     &\hspace{-0.5cm}=\! \Pi_{Z}\!\left[\!\left(\!\!\!\begin{array}{c} x(k) \! - \! \gamma_k\left(f_x\big(x(k)\big) \! + \! g_x\big(x(k)\big)^T\mu(k) \! + \! \alpha_kx(k)\right) \\ \mu(k) + \gamma_k\left(g\big(x(k)\big) - \alpha_k\mu(k)\right) \end{array}\!\!\!\right)\!\right] \label{eq:fullupdate}
\end{align}
until a fixed point $\hat{z}$ is reached. \hfill $\lozenge$
\end{algorithm}

Here $\alpha_k$ is the regularization parameter at timestep $k$ and $\gamma_k$ is the step-size
at the same timestep. 
In Section~\ref{sec:privacy} we will use Algorithm~\ref{alg:tikh} to solve a private
optimization problem,
and in Section~\ref{sec:conv} we provide hypotheses on $\gamma_k$ and $\alpha_k$
sufficient for convergence. Currently we show the applicability
of this style of solution to the cloud architecture mentioned above. 

\subsection{Communications} \label{ss:comms}
If we separate the update law in Algorithm~\ref{alg:tikh} to examine the per-agent (primal) update law, we find that agent $i$ executes
\begin{equation} \label{eq:per-agent-update}
x_i(k+1) = \Pi_{X_i}\Bigg[x_i(k) - \gamma_k\left(\nabla_i f_i(x_i(k)) + g_{x_i}(x(k))^T\mu(k) + \alpha_kx_i(k)\right)\Bigg].
\end{equation}
The only terms on the right-hand side of this update law that contain information from other agents are $g_{x_i}\big(x(k)\big)$ and $\mu(k)$. 
Though $g_{x_i}$ is a function of all states in the network, 
the agents do not send their states to each other directly to allow for its computation because doing so may reveal sensitive information.
Instead, every agent sends its state
to a trusted cloud computer which computes $g_{x_i}(x(k))$ for every $i \in I$. Because no agent
has every agent's state value, no agent can compute $\mu(k)$ (cf. Equation \eqref{eq:fullupdate}) and therefore the cloud computes $\mu(k)$ as well
using the update law
\begin{equation}
\mu(k+1) = \Pi_{\mathbb{M}}\left[\mu(k) + \gamma_k\left(g\big(x(k)\big) - \alpha_k\mu(k)\right)\right]. 
\end{equation}
Then, to use Algorithm~\ref{alg:tikh} with this architecture, the cloud
sends (private forms of) $g_{x_i}\big(x(k)\big)$ and $\mu(k)$ to agent $i$; the modifications to these quantities to make
them private are covered in Section~\ref{sec:privacy}. The cloud is assumed to be a powerful computer capable of carrying out
these calculations quickly so that they reliably arrive at the agents in a timely fashion. 

With this communications scheme, at timestep $k$ four actions occur. First, agent $i$ sends $x_i(k)$ to the cloud
and the cloud assembles all agents' states into the vector $x(k)$. Second, the cloud computes $\mu(k)$ and $g_{x_i}\big(x(k)\big)$ 
for all $i \in I$ in a differentially private way. 
Third, the cloud sends a private form of $g_{x_i}\big(x(k)\big)^T\mu(k)$
to agent $i$. Fourth, agent $i$ computes $x_i(k+1)$ while the cloud
simultaneously computes $\mu(k+1)$, and then this sequence of communications and computations
is repeated. Because this happens at every timestep, information in the network is always synchronized when computations
occur and there is no disagreement between the agents or cloud as to what the value of a particular state is. As a result, the 
computations that are spread across the network in this manner produce identical results to Algorithm~\ref{alg:tikh},
and the ensemble problem is, mathematically, equivalent to the cloud-based
multi-agent problem. 

For simplicity, the forthcoming analysis will be carried out in the ensemble setting. Despite the mathematical equivalence between the multi-agent
and ensemble approaches,
the advantage of the cloud-based approach in practice is that it allows for each agent's \new{state trajectory} to be kept private 
while the ensemble approach does not. 

\section{Private Optimization} \label{sec:privacy}
Differential privacy originates in the database literature
in computer science and was originally designed to keep
individual entries of a database private \cite{dwork13}. It
has recently been extended to the setting of dynamical systems in
\cite{pappas14}. Differential privacy offers a formal definition of privacy
as well as resilience to post-processing and robustness to side information.
This resilience to post-processing prevents an adversary from weakening the guarantees of 
differential privacy by performing post-hoc calculations on private
information. Robustness to side information guarantees that an adversary
cannot use information it has gleaned from an alternate source to 
fully defeat differential privacy. 
Below we first review differential privacy,
\new{
then give a formal private optimization problem statement,
} 
and finally discuss applying privacy to Problem~\ref{prob:vi}. 

\subsection{\new{Differentially Private Systems}}
Let there be $N$ input signals to a system, each contributed by some user.
The $i^{th}$ input signal is denoted $u_i$ and is contained in the 
set $\tl^{s_i}_{p_i}$, namely the space of sequences of $s_i$-vectors 
equipped with the $p_i$ norm, with $s_i, p_i \in \N$, 
such that every finite truncation of $u_i$ is in $\ell^{s_i}_{p_i}$. More explicitly,
let $u_{i}(k)$ denote the $k^{th}$ element of $u_{i}$ and define
\begin{equation}
P_{T}u_i = \begin{cases} u_i(k) &\mbox{for } k \leq T \\
                            0 &\mbox{otherwise}. \end{cases}
\end{equation}
Then we say $u_i \in \tl^{s_i}_{p_i}$ if and only if $P_{T}u_i$ has finite $p_i$-norm
for all values of $T$. 
Using this definition, the full input space to the system is
\begin{equation}
\tl^{s}_{p} = \prod_{i=1}^{N} \tl^{s_i}_{p_i},
\end{equation}
where the product is meant in the Cartesian sense, and the system produces 
outputs in
\begin{equation}
y \in \tilde{\ell}^{r}_{q}.
\end{equation}
In this paper we consider the cases where $p_i = 1$ for all $i \in I$ or $p_i = 2$
for all $i \in I$. In the case of $p_i = 1$, the full input space to the system
is $\tl^{s}_{1}$ and we use the ordinary $1$-norm on this space. For $p_i = 2$,
we likewise use the ordinary $2$-norm on $\tl^{s}_{2}$. 
While each of $\|\cdot\|_1$ and $\|\cdot\|_2$ will be used for both the
$1$-norm and $2$-norm on $\R^n$ and $\tl^s_p$, the intent of each symbol
can be discerned from its argument each time it is used. 

To implement differential privacy, we must specify which inputs we wish
to generate ``similar'' outputs. To do this, fix a real number $B > 0$
and define the binary symmetric adjacency relation 
$\Adj_B : \tl^s_p \times \tl^s_p \to \{0, 1\}$ as  
\begin{equation} \label{eq:adj}
\Adj_B(u, \tilde{u}) = 1 \Leftrightarrow \|u - \tilde{u}\|_p \leq B. 
\end{equation}
Two inputs $u$ and $\tilde{u}$ for which $\Adj(u, \tilde{u}) = 1$
are called ``adjacent.'' 

Towards making precise the notion of ``similar'' outputs, fix a probability space
$(\Omega, \mathcal{F}, \P)$ and let $\mathscr{B}^d$ denote the Borel
$\sigma$-algebra on $\R^d$. 
\new{
Differential privacy is enforced by a \emph{mechanism}, which is a map $M$
taking the form
\begin{equation}
M : \tl^{s}_{p} \times \Omega \to \tl^{r}_{q},
\end{equation}
and the role of a mechanism is to approximate a system whose inputs are sensitive
information. 
We now state the definition of a differentially private mechanism. In this definition,
we use a $\sigma$-algebra over $\lqr$, denoted $\Sigma_{q, r}$.\footnote{\new{An explicit construction
of this $\sigma$-algebra can be found in \cite[Section III-A]{pappas14}, though we avoid a lengthy
exposition on $\Sigma_{q, r}$ due to the relatively minor role its technical details play
in the current work.}}
\begin{definition} \label{def:fulleddef}
A mechanism $M : \lsp \times \Omega \to \lqr$ is $(\epsilon, \delta)$-differentially
private if and only if, for all adjacent $u, \tilde{u} \in \lsp$ we have
\begin{equation} \label{eq:fulleddef}
\mathbb{P}(M(u) \in S) \leq e^{\epsilon} \mathbb{P}(M(\tilde{u}) \in S) + \delta
\end{equation}
for all $S \in \Sigma_{q, r}$. \hfill $\lozenge$
\end{definition}
}

In Equation \eqref{eq:fulleddef} it is $\epsilon$ and $\delta$ that determine the privacy policy and 
smaller values of each imply a greater level of privacy for users. 
In general $\epsilon$ should be kept small and typical values for $\epsilon$
range from $0.1$ to $\ln 3$. On the other hand, $\delta$ should be kept as small as possible
because it allows for zero probability events for $M(\tilde{u})$ to have non-zero probability
for $M(u)$ and therefore can allow for important losses in privacy by making it easy for an
adversary to distinguish between outputs. 
Common values for $\delta$ range from $0$ to $0.05$; $(\epsilon, 0)$-differential privacy
is called $\epsilon$-differential privacy and, in general, $\epsilon$-differential privacy
is stronger than $(\epsilon, \delta)$-differential privacy precisely because of the aforementioned
losses in privacy that can come from $\delta > 0$. 
For this reason, $(\epsilon,\delta)$-differential privacy can be regarded as 
a $\delta$-approximate form of $\epsilon$-differential privacy \cite{dwork06b}. 
For a fixed value of $\epsilon$, the benefit of
using even small values of $\delta > 0$ is that the variance
of noise added can be reduced while maintaining ``almost'' the same level of privacy. 

\new{
\subsection{Private Optimization Problem Statement}
In the setting of Problem~\ref{prob:vi}, we want to protect the
state trajectory, $\bx = \big(x(k)\big)_{k \in \N}$, which is
a sensitive signal in $\lsp$, and in so doing we protect each individual agent's 
state trajectory; for agent $i$ this is $\bx_i \in \tl^{s_i}_{p_i}$. 
As discussed in Section~\ref{sec:problem}, keeping individual agents' state trajectories private
is necessary when the cloud computes $g_{x_i}$ and $\mu$ at each time $k$. To implement privacy
in these computations, we regard each $g_{x_i}$ as a deterministic, causal, memory-less dynamical system and
seek to make each such system differentially private. 
Similarly, we regard $g$ as a deterministic, causal, memory-less dynamical
system and seek to make it differentially private as well. Due to the post-processing property of differential
privacy, computing $\mu$ using a private form of $g$ also implies that $\mu$ keeps each $\bx_i$ private. 

As discussed in Section~\ref{sec:problem}, the agents do not communicate with each other at all and, instead,
each agent only sends its state to the cloud. The cloud handles all required centralized computations
and sends (privatized forms of) their results to the agents. 
Denoting by $\hat{g}_{x_i}$ and $\hat{g}$ the private forms of $g_{x_i}$ and $g$, respectively, 
at time $k$ the cloud sends to agent $i$ the vector
\begin{equation*}
\hat{p}_i(k) = \hat{g}_{x_i}\big(x(k)\big)^T\hat{g}\big(x(k)\big). 
\end{equation*}
We are interested in having a team of agents optimize by having the cloud
send agent $i$ only $\hat{p}_i(k)$ at time $k$. We require that $\hat{p}_i(k)$ protect
$\bx_i$ for all $i \in I$, and we implement privacy by approximating
$g_{x_i}$ (for all $i \in I$) and $g$ by differentially private mechanisms. 
Using this method of communications, we state the following problem that incorporates 
both optimization and privacy objectives, and respects the fact that the objectives
and constraints in this problem are sensitive data.

\begin{customprob}{1} \label{prob:withprivacy} (\emph{Private optimization})
Solve Problem~\ref{prob:vi} using Algorithm~\ref{alg:tikh} while
\begin{enumerate}[i.]
\item the agents communicate only with the cloud (i.e., there is no inter-agent communication)
\item the cloud makes the systems $g$ and $g_{x_i}$, $i \in I$
(whose inputs are the agents' state trajectories) differentially private
in the sense of Definition~\ref{def:fulleddef}
\item agent $i$ does not share $f_i$ or $X_i$ with any other agent or the cloud
\item the cloud does not share $g$ with any agent. \hfill $\lozenge$
\end{enumerate} 
\end{customprob}
}

Towards solving Problem~\ref{prob:withprivacy}, we now review mechanisms which implement differential privacy
for dynamical systems. 

%

\subsection{Privacy-Preserving Mechanisms}
To define a mechanism for enforcing differential privacy, we must first also define the sensitivity of
a system, which is used to determine the variance of noise that must be added in 
a privacy-preserving mechanism.
Letting $\mathcal{G}$ be a deterministic causal system,
the sensitivity of $\mathcal{G}$ is an upper bound on the distance
between $\mathcal{G}(u)$ and $\mathcal{G}(\tilde{u})$ whenever $\Adj_B(u, \tilde{u}) = 1$
holds. Formally we define the $\ell_p$ sensitivity of $\mathcal{G}$, denoted
$\Delta_p\mathcal{G}$, as
\begin{equation}
\Delta_{p}\mathcal{G} := \sup_{u,\tilde{u} : \Adj_B(u, \tilde{u}) = 1} \|\mathcal{G}(u) - \mathcal{G}(\tilde{u})\|_{p}. 
\end{equation}

The mechanism we will use for $\epsilon$-differential privacy is the Laplace mechanism, which adds noise drawn
from a Laplace distribution. Below we use the notation $\Lap(\mu, b)$ to denote the Laplace distribution with 
mean $\mu$ and scale parameter $b$. 
\begin{theorem} \label{thm:lap} (\cite[Theorem 4]{pappas14})
Let the adjacency relation defined in Equation \eqref{eq:adj} be used with $p = 1$ and let $\mathcal{G}$ be a system
with sensitivity $\Delta_{1}\mathcal{G}$. Let a constant $\epsilon \geq 0$ be given and recall that $r$ 
is the dimension of the output space. Then the mechanism
\begin{equation}
M(x) = \mathcal{G}(x) + w
\end{equation}
where $w(k) \sim \Lap(0, b/\epsilon)^r$ and $b \geq \Delta_{1}\mathcal{G}$ is $\epsilon$-differentially private.
\hfill $\blacksquare$
\end{theorem}


For $(\epsilon,\delta)$-differential privacy, we will use the Gaussian mechanism. 
Its definition requires that we first define $\kappa(\delta, \epsilon)$ using the $\mathcal{Q}$-function, 
\begin{equation}
\mathcal{Q}(y) := \frac{1}{\sqrt{2\pi}} \int_{y}^{\infty} e^{-\frac{v^2}{2}} dv.
\end{equation}
The function $\kappa(\delta, \epsilon)$ is defined  
for $\epsilon \geq 0$ and $0 < \delta < \frac{1}{2}$ as
\begin{equation}
\kappa(\delta, \epsilon) := \frac{1}{2\epsilon}\big(K_{\delta} + \sqrt{K_{\delta}^2 + 2\epsilon}\big),
\end{equation}
where $K_{\delta} = \mathcal{Q}^{-1}(\delta)$. We now define the Gaussian mechanism. 
\begin{theorem} \label{thm:gauss} (\cite[Theorem 3]{pappas14})
Let the adjacency relation defined in Equation \eqref{eq:adj} be used with $p = 2$ and let $\mathcal{G}$ be a system
with sensitivity $\Delta_{2}\mathcal{G}$, with constants $\epsilon \geq 0$ and $0 < \delta < \frac{1}{2}$ given and $r$ 
the dimension of the output space. Then the mechanism
\begin{equation}
M(x) = \mathcal{G}(x) + w
\end{equation}
where $w(k) \sim \mathcal{N}(0, \sigma^2I_r)$ is $(\epsilon, \delta)$-differentially private for
$\sigma \geq \kappa(\delta,\epsilon)\Delta_2\mathcal{G}$. \hfill $\blacksquare$
\end{theorem}

Theorems~\ref{thm:lap} and \ref{thm:gauss} provide a lower bound on the variance of each noise that is added and
we assume that these variances are also chosen to be finite. We now compute the sensitivities that are needed
to implement differential privacy in Problem~\ref{prob:withprivacy}.


\subsection{\new{Computing Sensitivities}}
\new{
In Problem~\ref{prob:withprivacy} it is desired to protect the value of $\bx$, including from agents in the network.
In the per-agent update law in Equation \eqref{eq:per-agent-update}, $x(k)$
appears in $g_{x_i}$ and $g_{x_i}$ must therefore be made private before the cloud sends $g_{x_i}(x(k))$ to agent $i$. 
To protect $\bx$ in this way, the cloud adds noise directly to $g_{x_i}\big(x(k)\big)$, and
the variance of noise that must be added depends on the sensitivity of $g_{x_i}$. 
To compute the sensitivity
of $g_{x_i}$ we regard it as a memoryless dynamical system and generalize it to act on entire signals of states. 
Recalling that $x(k) \in X$ for all $k$, under Assumption~\ref{as:set}, $X$ is bounded and therefore $x(k)$
is as well for all $k \in \N$. Then $\bx \in \tl^n_p$. 
}

We now overload the notation $g_{x_i}$ by allowing it to act on elements of $\tl^n_p$.
In particular, $g_{x_i}$ acts on elements of $\R^n$ as before and for \new{state trajectories} $\bx \in \tl^n_p$ we define
\begin{equation}
g_{x_i}(\bx) := \big(g_{x_i}(x(k))\big)_{k \in \N}. 
\end{equation}
We now fix a real scalar $B > 0$. 
For two \new{state trajectories}, $\bx, \tilde{\bx} \in \tl^n_p$ such that $\Adj_B(\bx, \tilde{\bx}) = 1$ holds, we 
compute the sensitivity of $g_{x_i}$ according to
\begin{multline}
\Delta_pg_{x_i} = \sup_{\bx, \tilde{\bx} : \Adj_B(\bx, \tilde{\bx}) = 1} \|g_{x_i}(\bx) - g_{x_i}(\tilde{\bx})\|_p \\ 
\leq \sup_{\bx, \tilde{\bx} : \Adj_B(\bx, \tilde{\bx}) = 1} K^i_p \sqrt[p]{\sum_{k=0}^{\infty} \|x(k) - \tilde{x}(k)\|_p^p} \leq K^i_pB,
\end{multline}
where we have used $\|\bx - \tilde{\bx}\|_p \leq B$ and where this bound on the sensitivity holds for $g_{x_i}$ for all $i \in I$. 

In computing $\mu(k)$, the cloud must also add noise in some fashion because $\mu(k)$ depends
upon $x(k)$. We regard $g$ 
as a dynamical system and make it private, and the resilience of differential privacy to post-processing guarantees
that $\bmu = \big(\mu(k)\big)_{k \in \N}$ keeps $\bx$ private.
To compute the sensitivity of $g$ we extend it to act on $\bx \in \tl^n_p$ as above. For $\bx, \tilde{\bx} \in \tl^n_p$ satisfying
$\Adj_B(\bx, \tilde{\bx}) = 1$ we use the same procedure as was used above for $g_{x_i}$ to find
\begin{align}
\Delta_{p}g &= \sup_{\bx, \tilde{\bx} : \Adj_B(\bx, \tilde{\bx}) = 1} \|g(\bx) - g(\tilde{\bx})\|_{p} \leq K^g_pB.
\end{align}

\new{Having computed the requisite sensitivities, we return to solving Problem~\ref{prob:withprivacy}.} 

\subsection{Optimizing in the Presence of Noise} \label{ss:opt-noise}
We now examine how noise appears in Algorithm~\ref{alg:tikh} once it has
been added for privacy. For $g_{x_i}\big(x(k)\big)$ we add noise $w_i(k) \in \R^{m \times n_i}$ drawn from either a Laplace
or Gaussian distribution and for $g\big(x(k)\big)$ we add noise $w_g(k) \in \R^m$ drawn from the same class of distribution as the $w_i$, with
all noises independent. Define $w_x$ by
\begin{equation}
w_x = (w_1 \,\, w_2 \,\, \cdots \,\, w_n) \in \R^{m \times n}. 
\end{equation}

In ensemble form the private dynamics under consideration are
\begin{align}
z(k+1) &= \left(\begin{array}{c} x(k+1) \\ \mu(k+1) \end{array}\right) \\ &= \Pi_{Z}\left[z(k) - \gamma_k \left(\begin{array}{c} f_x\big(x(k)\big) + \left(\frac{\partial g}{\partial x}\big(x(k)\big) + w_x(k)\right)^T\mu(k) + \alpha_kx(k)\\ -g\big(x(k)\big) + w_g(k) + \alpha_k\mu(k)\end{array}\right)\right].
\end{align}
Expanding, we find
\begin{align} \label{eq:noisyupdate}
z(k+1) &= \left(\begin{array}{c} x(k+1) \\ \mu(k+1) \end{array}\right) \\ &= \Pi_{Z}\left[z(k) 
 - \gamma_k \left( 
\begin{array}{c} f_x\big(x(k)\big)  +  \frac{\partial g}{\partial x}\big(x(k)\big)^T\mu(k) + w_x(k)^T\mu(k) + \alpha_kx(k)  \\ -g\big(x(k)\big) + w_g(k) + \alpha_k\mu(k)\end{array}\right)\right].
\end{align}
Because $\mu(k) \in \R^m_{+}$, each element of $w_x(k)^T\mu(k)$ is some weighted combination of 
elements of $w_x(k)$ with non-negative weights. Combined with the independence of the noises
used for privacy, this results in each entry of $w_x(k)^T\mu(k)$ being a random variable having
variance that is the weighted sum of variances of elements of $w_x(k)$. With this in mind we
define the random vector $w_s(k) = w_x(k)^T\mu(k)$ (which we note has finite variance since $\mu(k)$ is contained in $\mathbb{M}$ and $w_x(k)$ has finite variance),
and zero mean (because $w_i(k)$ has zero mean for all $k \in \N$ and all $i \in I$). 
Then we can 
rewrite Equation \eqref{eq:noisyupdate} as
\begin{align}
z(k+1) &=  \Pi_{Z}\left[z(k) 
- \gamma_k \left(\begin{array}{c} f_x\big(x(k)\big) + \frac{\partial g}{\partial x}\big(x(k)\big)^T\mu(k) + w_s(k) + \alpha_kx(k)
 \\ -g\big(x(k)\big) + w_g(k) + \alpha_k\mu(k)\end{array}\right)\right] \\
    &= \Pi_{Z}\left[z(k) - \gamma_k\left(G\big(z(k)\big) + \alpha_kz(k) + w(k)\right)\right],
\end{align}
where $w(k)$ denotes the noise added at timestep $k$ and aggregates all noisy signals used for privacy.
We state this stochastic update law as Algorithm~\ref{alg:stoch}.
\begin{algorithm} \label{alg:stoch}
Given $z(0) \in Z$, apply the update law
\begin{equation}
z(k+1) = \Pi_{Z}\left[z(k) - \gamma_k\left(G\big(z(k)\big) + \alpha_kz(k) + w(k)\right)\right]
\end{equation}
until a fixed point $\hat{z} \in Z$ is reached. \hfill $\lozenge$
\end{algorithm}

We note that by its definition $\E[w(k)] = 0$, and 
observe that this noisy update law is equivalent to Algorithm~\ref{alg:tikh} with an additional noise term added.

\section{Convergence of Private Optimization} \label{sec:conv}
In this section we prove the convergence of Algorithm~\ref{alg:stoch}. 
Algorithm~\ref{alg:stoch} was first presented in \cite{bakushinskii74} without noise
and was presented in its noisy form in \cite{poljak78}. Both papers omit proofs and, due to the
heavy dependence of this work upon Algorithm~\ref{alg:stoch} we provide a proof here. To the best
of our knowledge a proof of the convergence of Algorithm~\ref{alg:stoch} as stated in \cite{poljak78}
is not available in the literature; similar work is presented in \cite{yousefian13,koshal13} which cover
algorithms related to Algorithm~\ref{alg:stoch}, though those works impose additional assumptions
upon $\alpha_k$ and $\gamma_k$ due to the differences in the problems studied in those works. 

\subsection{Main Convergence Result}
Now we explore in depth solving variational inequalities using a Tikhonov regularized projection method, the basic
elements of which are covered in \cite[Section 12.2]{facchinei03}. 
Earlier it was stated that if $SOL(K, F) \neq \0$, then for $\xi_k \in SOL(K, F + \alpha_kI)$
we have $\xi_k \to z_0$ where $z_0$ is the least-norm element of $SOL(K, F)$. Using that
$\{\xi_k\}_{k \in \N}$ is a convergent sequence, we find that 
$\{\|\xi_k\|\}_{k \in \N}$ is bounded and, in particular, there is some $M_{\xi}$
such that $\|\xi_k\| \leq M_{\xi}$ for all $k$, e.g., $\|\xi_k\| \leq \sup_{z \in Z} \|z\|$. 

Using this fact, the following lemma relates points $z(k)$ generated by Algorithm~\ref{alg:stoch} 
to successive solutions to the problems $VI(Z, G + \alpha_kI)$ (each with $\alpha_k$ held constant).
Recalling that $\xi_k$ is the unique solution to $VI(Z, G + \alpha_kI)$, we have the following result. 
\begin{lemma} \label{lem:bigineq}
For all $k \in \N$
\begin{equation}
\|z(k) - \xi_k\|^2 \leq (1 + \gamma_k\alpha_k)\|z(k) - \xi_{k-1}\|^2 
 + M_{\xi}^2\left(\frac{\alpha_{k-1} - \alpha_k}{\alpha_k}\right)^2\left(\frac{1 + \gamma_k\alpha_k}{\gamma_k\alpha_k}\right).
\end{equation}
\end{lemma}
\emph{Proof:}  First note
that because $\xi_k$ solves $VI(K, G + \alpha_k I)$, we have
\begin{equation} \label{eq:xi1}
(\xi_{k-1} - \xi_{k})^T\big(G(\xi_k) + \alpha_k\xi_k\big) \geq 0. 
\end{equation}
Similarly for $\xi_{k-1}$ we find
\begin{equation} \label{eq:xi2}
(\xi_k - \xi_{k-1})^T\big(G(\xi_{k-1}) + \alpha_{k-1}\xi_{k-1}\big) \geq 0.
\end{equation}
Summing Equations \eqref{eq:xi1} and \eqref{eq:xi2}, and using the monotonicity of $G$ gives
\begin{equation}
(\xi_{k-1} - \xi_k)^T(\alpha_k \xi_k - \alpha_{k-1}\xi_{k-1}) \geq 0. 
\end{equation}

Adding and subtracting $\alpha_k\xi_{k-1}$ inside the second set of parentheses then gives
\begin{multline}
(\xi_{k-1} - \xi_k)^T(\alpha_k I - \alpha_{k-1}I)\xi_{k-1} \\ \geq \alpha_k (\xi_{k-1} - \xi_k)^T(\xi_{k-1} - \xi_k)
 = \alpha_k \|\xi_{k-1} - \xi_k\|^2.
\end{multline}
Using the Cauchy-Schwarz inequality results in
\begin{equation} \label{eq:xidiff}
\|\xi_{k-1} - \xi_k\| \leq \frac{|\alpha_{k-1} - \alpha_k|}{\alpha_k} M_{\xi}.
\end{equation}

Expanding the term $\|z(k) + \xi_{k-1} - \xi_{k-1} - \xi_k\|^2$ and applying Equation \eqref{eq:xidiff} then gives
\begin{multline} \label{eq:lemmaalmost}
\|z(k) - \xi_k\|^2 \leq\|z(k) - \xi_{k-1}\|^2+\left(\frac{\alpha_{k-1} - \alpha_k}{\alpha_k}\right)^2 M_{\xi}^2
\\ + 2\left(\frac{|\alpha_{k-1} - \alpha_k|}{\alpha_k}\right)M_{\xi}\|z(k) - \xi_{k-1}\|.
\end{multline}
For the third term on the right-hand side above we have
\begin{multline}
2\left(\frac{|\alpha_{k-1} - \alpha_k|}{\alpha_k}\right)M_{\xi}\|z(k) - \xi_{k-1}\|  \\ 
= 2\sqrt{\gamma_k\alpha_k}\|z(k) 
- \xi_{k-1}\|\left(\frac{|\alpha_{k-1} - \alpha_k|}{\alpha_k \sqrt{\gamma_k\alpha_k}}\right)M_{\xi} \\
\leq \gamma_k\alpha_k\|z(k) - \xi_{k-1}\|^2 + M_{\xi}^2\left(\frac{(\alpha_{k-1} - \alpha_k)^2}{\gamma_k\alpha_k^3}\right)\label{eq:bigineqend},
\end{multline} 
where we have used that fact that $a^2 + b^2 \geq 2ab$ for $a, b \in \R$.

Substituting Equation \eqref{eq:bigineqend} into Equation \eqref{eq:lemmaalmost} gives the desired result. \hfill $\blacksquare$

The other lemma we need concerns the convergence of sequences of random variables and enables a 
Lyapunov-like argument to be made for their convergence. 
\begin{lemma} \label{lem:lyap1} (\cite{poljak87}, Lemma $10$, Page $49$)
Let $v_0, \ldots, v_k$ be a sequence of independent random variables with $v_k \geq 0$ and $\E[v_0] < \infty$. Suppose that
\begin{equation}
\E[v_{k+1}] \leq (1 - \tau_k)v_k + \sigma_k,
\end{equation}
with 
\begin{equation}
0 \leq \tau_k \leq 1, \quad \sigma_k \geq 0, \quad \sum_{k=0}^{\infty} \tau_k = \infty,  \quad \frac{\sigma_k}{\tau_k} \to 0.
\end{equation}
Then $\E[v_k] \to 0$. If, in addition, we have
\begin{equation}
\sum_{k=0}^{\infty} \sigma_k < \infty,
\end{equation}
then $v_k \to 0$ almost surely and 
\begin{equation}
\P(v_j \leq \epsilon \textnormal{ for all } j \geq k) \geq 1 - \frac{1}{\epsilon}\left(\E[v_k] + \sum_{i=k}^{\infty} \sigma _i\right).
\end{equation} \hfill $\blacksquare$
\end{lemma}



We now prove the convergence of Algorithm~\ref{alg:stoch}. 

\begin{theorem} \label{thm:conv}
Let Assumptions 1-4 hold. Suppose that ${\gamma_k > 0}$ and $\alpha_k > 0$ satisfy the following four conditions: 
\begin{enumerate}
\item $\sum_{k=0}^{\infty} \gamma_k\alpha_k = \infty$ 
\item $\frac{\gamma_k}{\alpha_k} \to 0$ 
\item $\alpha_k \to 0$
\item $\frac{(\alpha_{k-1} - \alpha_k)}{\gamma_k\alpha_k^2} \to 0$.
\end{enumerate}
Then for noise signal $w$ with $\E[w(k)] = 0$ and bounded variance
for all $k \in \N$, for the update rule
\begin{equation}
z(k+1) = \Pi_{Z}\left[z(k) - \gamma_k\left(G\big(z(k)\big) + \alpha_kz(k) + w(k)\right)\right]
\end{equation}
we have $\E[\|z(k) - z_0\|^2] \to 0$, where $z_0$ is the least-norm solution to Problem~\ref{prob:vi}. 

Let $L_G$ be the Lipschitz constant of $G$. If, in addition
to the above, the sequence of terms
\begin{multline}
\sigma_k := \left(1 - \gamma_k\alpha_k\left(2 - \gamma_k\alpha_k - \frac{\gamma_k}{\alpha_k}L_G^2 - 2\gamma_kL_G\right)\right)M_{\xi}^2
\cdot \left(\frac{\alpha_{k-1} - \alpha_k}{\alpha_k}\right)^2\left(\frac{1+\gamma_k\alpha_k}{\gamma_k\alpha_k}\right) \\ + \gamma_k^2\E[\|w(k)\|^2]
\end{multline}
is summable, then the convergence estimate 
\begin{equation} \label{eq:rateestim}
\P(\|z(j) - \xi_{j-1}\|^2 \leq \epsilon \textnormal{ for all } j \geq k) \geq 1 - \frac{1}{\epsilon}\left(\E[\|z(k) - \xi_{k-1}\|^2] + \sum_{i=k}^{\infty} \sigma _i\right)
\end{equation}
holds.
\end{theorem}
\emph{Proof:} It was established in Section~\ref{ss:problem} that $SOL(Z, G) \neq \emptyset$ so that
$\xi_k \to z_0$ where $z_0$ is the least-norm element of $SOL(Z, G)$ and where $\xi_k$ solves $VI(Z, G + \alpha_kI)$.
We now show that $z(k+1) \to \xi_k$. 

Because $\xi_k$ solves $VI(Z, G + \alpha_kI)$ we have
\begin{equation}
\xi_k = \Pi_{Z}\left[\xi_k - \gamma_k\big(G(\xi_k) + \alpha_k\xi_k\big)\right].
\end{equation}
Using the non-expansive property of the projection operator and taking the expectation of both sides we find
\begin{multline}
\E\big[\|z(k+1) - \xi_k\|^2\big] = \\
\E\bigg[\Big\|\Pi_{Z}\Big[z(k) - \gamma_k\big(G(z(k)) + \alpha_kz(k) + w(k)\big)\Big] 
- \Pi_{Z}\Big[\xi_k - \gamma_k\big(G(\xi_k) + \alpha_k\xi_k\big)\Big]\Big\|^2\bigg] \\
 \leq \E\bigg[\|z(k) - \xi_k + \gamma_k(G(\xi_k) - G(z(k))) - \gamma_k\alpha_k(z(k) - \xi_k) - \gamma_kw(k)\|^2 \bigg] \\
 \leq \|z(k) - \xi_k\|^2 - 2\gamma_k\alpha_k\|z(k) - \xi_k\|^2 + \gamma_k^2\|G(\xi_k) - G(z(k))\|^2 \\
 + 2\gamma_k^2\alpha_k(G(\xi_k) - G(z(k)))^T(\xi_k - z(k))  + \gamma_k^2\alpha_k^2\|\xi_k - z(k)\|^2 + \gamma_k^2\E[\|w(k)\|^2]
\end{multline}
where the last inequality follows from the monotonicity of $G$, and where the fact that $\E[w(k)] = 0$ has caused all terms containing $w(k)$
except $\E[\|w(k)\|^2]$ to vanish. 

Using the Cauchy-Schwarz inequality then gives
\begin{multline}
\E[\|z(k+1) - \xi_k\|^2] \leq \|z(k) - \xi_k\|^2 - 2\gamma_k\alpha_k\|z(k) - \xi_k\|^2 
         \\ + \gamma_k^2\|G(\xi_k) - G(z(k))\|^2 + 2\gamma_k^2\alpha_k\|G(\xi_k) - G(z(k))\| \|\xi_k - z(k)\| 
         \\ + \gamma_k^2\alpha_k^2\|\xi_k - z(k)\|^2 + \gamma_k^2\E[\|w(k)\|^2].
\end{multline}

Assumptions~\ref{as:objectives}-\ref{as:constraints} and the compactness
of $\mathbb{M}$ together imply that $G$ is Lipschitz and, denoting its Lipschitz
constant by $L_G$, we have
\begin{multline} \label{eq:bigz1}
\E[\|z(k+1) - \xi_k\|^2] \leq \\ \Big(1 - 2\gamma_k\alpha_k + \gamma_k^2L_G^2 + 2\gamma_k^2\alpha_kL_G + \gamma_k^2\alpha_k^2\Big)\|z(k) - \xi_k\|^2 + \gamma_k^2\E[\|w(k)\|^2].
\end{multline}

Defining
\begin{equation}
\theta_k := 1 - \gamma_k\alpha_k\left(2 - \gamma_k\alpha_k - \frac{\gamma_k}{\alpha_k}L_G^2 - 2\gamma_kL_G\right)
\end{equation}
and
\begin{equation}
\rho_k := M_{\xi}^2\left(\frac{\alpha_{k-1} - \alpha_k}{\alpha_k}\right)^2\left(\frac{1+\gamma_k\alpha_k}{\gamma_k\alpha_k}\right),
\end{equation}
applying Lemma~\ref{lem:bigineq} then gives
\begin{equation} \label{eq:t1bigone}
\E[\|z(k+1) - \xi_k\|^2] \leq \theta_k(1 + \gamma_k\alpha_k)\|z(k) - \xi_{k-1}\|^2 + \theta_k\rho_k + \gamma_k^2\E[\|w(k)\|^2].
\end{equation}

By hypothesis we have
\begin{equation}
\gamma_k\alpha_k \to 0, \quad \frac{\gamma_k}{\alpha_k} \to 0, \quad \gamma_k \to 0,
\end{equation}
with $\alpha_k > 0$ and $\gamma_k > 0$ for all $k$. 
Then there exists an $M > 0$ such that for all $k \geq M$ we have
\begin{equation}
\gamma_k\alpha_k \in (0, 1) \textnormal{ and } 0 \leq 1 - 2\gamma_k\alpha_k \leq \theta_k \leq 1 - \gamma_k\alpha_k.
\end{equation}
Then for all $k \geq M$
\begin{equation}
\theta_k(1 + \gamma_k\alpha_k) \leq \theta_k + \gamma_k\alpha_k
\end{equation}
and thus for all $k \geq M$ 
\begin{equation}
\theta_k(1 + \gamma_k\alpha_k) \leq 1 - \gamma_k\alpha_k(1 - \gamma_k\alpha_k - \frac{\gamma_k}{\alpha_k}L_G^2 - 2\gamma_kL_G) \in (0,1).
\end{equation}

In particular, take some $\theta \in (0, 1)$ so that 
\begin{equation}
1 - \gamma_k\alpha_k\left(1 - \gamma_k\alpha_k - \frac{\gamma_k}{\alpha_k}L_G^2 - 2\gamma_kL_G\right) \leq 1 - \gamma_k\alpha_k\theta
\end{equation}
for all $k \geq M$. 
Then $1 - \gamma_k\alpha_k\theta \in (0, 1)$. Setting $\tau_k = \gamma_k\alpha_k\theta$ and
$\sigma_k = \rho_k\theta_k + \gamma_k^2\E[\|w(k)\|^2]$ we rewrite Equation \eqref{eq:t1bigone} as 
\begin{equation} \label{eq:recur}
\E[\|z(k+1) - \xi_k\|^2] \leq (1 - \tau_k)\|z(k) - \xi_{k-1}\|^2 + \sigma_k.
\end{equation}

All that remains is to show that the conditions of Lemma~\ref{lem:lyap1} are met. 
First, $\tau_k \in (0, 1)$ by construction. For all $k \geq M$ we have
$\rho_k \geq 0$ and $\theta_k \geq 0$ so that $\sigma_k \geq 0$. Regarding
summability of $\tau_k$ we find
\begin{align}
\sum_{k = M}^{\infty} \tau_k = \theta\sum_{k=M}^{\infty} \gamma_k\alpha_k = \infty
\end{align}
by hypothesis. 
To show that $\sigma_k/\tau_k \to 0$ we have
\begin{align}
\frac{\sigma_k}{\tau_k} &= \frac{\theta_k \left(\frac{\alpha_{k-1} - \alpha_k}{\alpha_k}\right)^2
\left(\frac{1 + \gamma_k\alpha_k}{\gamma_k\alpha_k}\right)M_{\xi}^2 + \gamma_k^2\E[\|w(k)\|^2]}{\gamma_k\alpha_k\theta} \\
&=  (1 - 2\gamma_k\alpha_k + \gamma_k^2\alpha_k^2 + \gamma_k^2L_G^2 + 2\alpha_k\gamma_k^2L_G)M_{\xi}^2 \\
&\cdot \left(\frac{\alpha_{k-1} - \alpha_k}{\gamma_k\alpha_k^2}\right)^2(1 + \gamma_k\alpha_k) \frac{1}{\theta} + \frac{1}{\theta}\frac{\gamma_k}{\alpha_k}\E[\|w(k)\|^2].
\label{eq:seqbigone}
\end{align}
Using the hypotheses regarding $\gamma_k$ and $\alpha_k$ we have
\begin{equation}
1 - 2\gamma_k\alpha_k + \gamma_k^2\alpha_k^2 + \gamma_k^2L_G^2 + 2\alpha_k\gamma_k^2L_G \to 1, \quad
1 + \gamma_k\alpha_k \to 1,
\end{equation}
along with
\begin{equation}
\left(\frac{\alpha_{k-1} - \alpha_k}{\gamma_k\alpha_k^2}\right)^2 \to 0,
\end{equation}
so that the first term in $\frac{\sigma_k}{\tau_k}$ goes to zero. It was established in Section~\ref{ss:opt-noise}
that $\E[\|w(k)\|^2]$ is bounded above for all $k$, namely that $\E[\|w(k)\|^2] \leq K_{w}$ for some $K_w > 0$.
Because $\frac{\gamma_k}{\alpha_k} \to 0$ we have $\frac{\gamma_k}{\alpha_k}K_w \to 0$ and hence
$\frac{\sigma_k}{\tau_k} \to 0$
as desired and the first part of the theorem follows from Lemma~\ref{lem:lyap1}. 

When the sequence $\{\sigma_k\}_{k \in \N}$ is summable, the additional
convergence estimate is a straightforward application of Lemma~\ref{lem:lyap1} as well. \hfill $\blacksquare$

One valid choice of $\gamma_k$ and $\alpha_k$ satisfying conditions 1-4 in Theorem~\ref{thm:conv} is
\begin{equation}
\alpha = \bar{\alpha}k^{-c_1} \textnormal{ and } \gamma = \bar{\gamma}k^{-c_2},
\end{equation}
with $\bar{\alpha} > 0$, $\bar{\gamma} > 0$, $0 < c_1 < c_2$, and $c_1 + c_2 < 1$ \cite{poljak78}. 

\subsection{Convergence Rate Estimates}
For the above choice of step-size, we derive bounds on $c_1$ and $c_2$ which are sufficient
to make $\sigma_k$ summable. As shown in the proof of Theorem~\ref{thm:conv}, 
there exists an $M > 0$ such that for all $k \geq M$ we have
$\theta_k \in (0, 1)$,
so that for all $k \geq M$ we have
\begin{equation} \label{eq:bound}
\sigma_k \leq M_{\xi}^2\left(\frac{\alpha_{k-1} - \alpha_k}{\alpha_k}\right)^2\left(\frac{1 + \gamma_k\alpha_k}{\gamma_k\alpha_k}\right) + \gamma_k^2 \E\left[\|w(k)\|^2\right].
\end{equation}
To make the second term in Equation \eqref{eq:bound} summable, we can set 
\begin{equation}
\gamma_k = \bar{\gamma}k^{-c_2}
\end{equation}
with $c_2 > \frac{1}{2}$. Again using $K_w$ to denote an upper bound on the variance of $w(k)$ gives
\begin{equation} \label{eq:sum1}
\sum_{k=1}^{\infty} \gamma_k^2\E[\|w(k)\|^2] \leq \bar{\gamma}^2K_w\zeta(2c_2),
\end{equation}
where $\zeta(\cdot)$ is the Riemann zeta function \cite{edwards01}, defined as
\begin{equation}
\zeta(p) = \sum_{n=1}^{\infty} \frac{1}{n^p},
\end{equation} 
which takes finite values for arguments $p > 1$. 

Regarding the first term in Equation \eqref{eq:bound}, we note that there is some $\hat{M} > 0$ such that 
\begin{equation} \label{eq:kless}
1 \leq \frac{1}{\bar{\alpha}\bar{\gamma}k^{-c_1}k^{-c_2}}
\end{equation}
for all $k \geq \hat{M}$
and therefore
\begin{equation}
\left(\frac{\alpha_{k-1} - \alpha_k}{\alpha_k}\right)^2\left(\frac{1 + \gamma_k\alpha_k}{\gamma_k\alpha_k}\right) \leq 
2\left(\frac{\alpha_{k-1} - \alpha_k}{\alpha_k}\right)^2\frac{1}{\bar{\alpha}\bar{\gamma}k^{-(c_1 + c_2)}}
\end{equation}
for all such $k$. Substituting $\alpha_k = \bar{\alpha}k^{-c_1}$ and
expanding the squared term gives
\begin{equation} \label{eq:needexpand}
2\left(\frac{\alpha_{k-1} - \alpha_k}{\alpha_k}\right)^2\frac{1}{\bar{\alpha}\bar{\gamma}k^{-(c_1 + c_2)}} = 
2\frac{\left[\left(1 - \frac{1}{k}\right)^{-2c_1} - 2\left(1 - \frac{1}{k}\right)^{-c_1} + 1\right]}{\bar{\alpha}\bar{\gamma}k^{-(c_1 + c_2)}}.
\end{equation}
To approximate the terms containing $1/k$ we use a (truncated) power series expansion, namely that for $x \in (-1, 1)$
\begin{equation} \label{eq:expansion}
(1 - x)^{-r} \approx 1 + rx + \frac{1}{2}r(r+1)x^2 + \frac{1}{6}r(r+1)(r+2)x^3.
\end{equation}
Applying Equation \eqref{eq:expansion} to Equation \eqref{eq:needexpand} gives
\begin{equation}
2\frac{\left[\left(1 - \frac{1}{k}\right)^{-2c_1} - 2\left(1 - \frac{1}{k}\right)^{-c_1} + 
1\right]}{\bar{\alpha}\bar{\gamma}k^{-(c_1 + c_2)}} \approx \frac{2c_1^2}{\bar{\alpha}\bar{\gamma}k^{2 - (c_1 + c_2)}} + \frac{2(c_1^3 + c_1^2)}{\bar{\alpha}\bar{\gamma}k^{3 - (c_1 + c_2)}}.
\end{equation}
We see that sums of such terms are given by
\begin{multline} \label{eq:sum2}
\sum_{k=1}^{\infty} 2\frac{\left[\left(1 - \frac{1}{k}\right)^{-2c_1} - 2\left(1 - \frac{1}{k}\right)^{-c_1} + 1\right]}{\bar{\alpha}\bar{\gamma}k^{-(c_1 + c_2)}} \\ \approx 
\frac{2c_1^2}{\bar{\alpha}\bar{\gamma}}\zeta(2 - (c_1 + c_2)) + \frac{2(c_1^3 + c_1^2)}{\bar{\alpha}\bar{\gamma}}\zeta(3 - (c_1 + c_2)).
\end{multline}
Returning to Equation \eqref{eq:bound} and using the results of Equations \eqref{eq:sum1} and \eqref{eq:sum2} gives
\begin{equation} \label{eq:sumsigma}
\sum_{k=1}^{\infty} \sigma_k \lessapprox \bar{\gamma}^2K_w\zeta(2c_2) + \frac{2c_1^2}{\bar{\alpha}\bar{\gamma}}\zeta(2 - (c_1 + c_2)) + \frac{2(c_1^3 + c_1^2)}{\bar{\alpha}\bar{\gamma}}\zeta(3 - (c_1 + c_2)).
\end{equation}

Due to the approximations made and ranges of $k$ considered in bounding this sum, we can only guarantee that
the convergence estimate relying on $\sum_{k=1}^{\infty} \sigma_k$ will hold for $k \geq \max\{M, \hat{M}\}$. 
However, for $k \leq M$ we will often have $\sigma_k < 0$ (as when $L_G$ is large) and thus we expect the bound in Equation 
\eqref{eq:sumsigma} to hold for a range of values of $k \leq M$ because negative terms with such indices have been over-estimated by including positive 
terms at such indices in Equation \eqref{eq:sumsigma}. In addition, we expect $\bar{\alpha}$ and $\bar{\gamma}$ to be small enough
that $\hat{M}$ will often be small, e.g., less than $10$, thus allowing this bound to hold over a wide range of values of $k$. 

\new{
To apply the bound in Equation \eqref{eq:rateestim}, we also need to estimate the term $\E[\|z(k) - \xi_{k-1}\|^2]$. 
Returning to Equation \eqref{eq:recur} and taking the expectation of both sides one timestep earlier gives
\begin{equation} \label{eq:exprecur}
\E[\|z(k) - \xi_{k-1}\|^2] \leq (1 - \tau_{k-1})\E[\|z(k-1) - \xi_{k-2}\|^2] + \sigma_{k-1},
\end{equation}
which is a (time-varying) affine recurrence relation in the expected error in the optimization algorithm. 
Solving Equation \eqref{eq:exprecur} (see e.g., \cite{greene81}, Section 2.1.1.2), we find that
\begin{equation}
\E[\|z(k) - \xi_{k-1}\|^2] = \left(\prod_{n=0}^{k-1} (1 - \tau_n)\right)\left(\E[\|z(1) - \xi_0\|^2] + \sum_{m=0}^{k-1} \frac{\sigma_m}{\prod_{i=0}^{m}(1 - \tau_i)}\right).
\end{equation}
Defining the diameter of the set $Z$ via $D_z = \sup_{z_1, z_2 \in Z} \|z_1 - z_2\|$, we can bound the initial error via
$\E[\|z(1) - \xi_0\|^2] \leq D_z^2$,
giving
\begin{equation} \label{eq:expbound}
\E[\|z(k) - \xi_{k-1}\|^2] \leq \left(\prod_{n=0}^{k-1} (1 - \tau_n)\right)\left(D_z^2 + \sum_{m=0}^{k-1} \frac{\sigma_m}{\prod_{i=0}^{m}(1 - \tau_i)}\right).
\end{equation}

In Equation \eqref{eq:rateestim}, one can compute the sum $\sum_{i=k}^{\infty} \sigma_i$ by using the analytic bound for 
$\sum_{k=1}^{\infty} \sigma_k$ given in Equation \eqref{eq:sumsigma} and subtracting the first $k-1$ values of $\sigma_i$
from this value. Combined with the upper bound on $\E[\|z(k) - \xi_{k-1}\|^2]$ given in Equation \eqref{eq:expbound}, one can then
use the bound in Equation \eqref{eq:rateestim} to determine the probability with which the error in the optimization algorithm stays
within some bound for all time. 
}
\new{Having explored convergence in the presence of privacy, we now examine the trade-off between the two competing objectives of privacy and convergence.}

\subsection{\new{The Trade-off Between Privacy and Convergence}}
\new{In this section we derive a quantifiable trade-off between privacy and convergence, and for concreteness
we focus on the case of $\epsilon$-differential privacy, though a similar trade-off can be derived for
$(\epsilon, \delta)$-differential privacy. 

Returning to Equation \eqref{eq:t1bigone} we find the inequality
\begin{equation} \label{eq:trade1}
\E[\|z(k+1) - \xi_k\|^2] \leq \theta_k(1 + \alpha_k\gamma_k)\|z(k) - \xi_{k-1}\|^2 + \theta_k\rho_k + \gamma_k^2\E[\|w(k)\|^2],
\end{equation}
where we see that only the term $\E[\|w(k)\|^2]$ depends upon the noise added for privacy. Given that
$w(k)$ has zero mean, we find $\E[\|w(k)\|^2] = \textnormal{var}\big(w(k)\big)$. In the case of
$\epsilon$-differential privacy we have $w(k) \sim \Lap(0, b/\epsilon)^r$ so that 
\begin{equation} \label{eq:trade2}
\textnormal{var}[w(k)] = \frac{W}{\epsilon^2},
\end{equation}
where $W := W(\Delta_1 g, \Delta_1 g_{x_i}, B)$ is a constant that
depends upon the systems of interest, $g$ and $g_{x_i}$, and the adjacency parameter, $B$. 

Returning to Equation \eqref{eq:trade1} and substituting in Equation \eqref{eq:trade2} we find
\begin{equation} \label{eq:trade4}
\E[\|z(k+1) - \xi_k\|^2] \leq \theta_k(1 + \alpha_k\gamma_k)\|z(k) - \xi_{k-1}\|^2 + \theta_k\rho_k + \gamma_k^2\frac{W}{\epsilon^2}.
\end{equation}
The additive term $\gamma_k^2\frac{W}{\epsilon^2}$ is the only term in which the privacy parameter
$\epsilon$ appears, and this term can be regarded as a penalty on convergence because it 
allows the expected error $\E[\|z(k+1) - \xi_k\|^2]$ to grow from $\|z(k) - \xi_{k-1}\|$. Viewing this term as a convergence
penalty then reveals a fundamental trade-off between privacy and convergence: implementing $\epsilon$-differential
privacy comes at the cost of a convergence penalty proportional to $1/\epsilon^2$. We state this trade-off succinctly and
informally by writing
\begin{equation}
\textnormal{Privacy}(\epsilon) \Longleftrightarrow \textnormal{Convergence}(\epsilon^{-2}). 
\end{equation}

}
\vspace{-0.4cm}

\section{Simulation Results} \label{sec:simulation}
Below we present numerical simulation results for a system
with $n = 10$ agents and $m = 6$ constraints. We simulate both
$\epsilon$- and $(\epsilon, \delta)$-differential privacy.

\vspace{-0.3cm}

\subsection{Example Problem}
Let there be $n = 10$ agents, each with state $x_i \in \R^2$
and using ensemble objective function
\vspace{-0.3cm}
\begin{multline}
f(x) = \big((x_{1, 1} - 5) + (x_{1, 2} + 5)\big) + \|x_2\|^2 + \left\|x_3 - \left(\begin{array}{r} -7 \\ 7 \end{array}\right) \right\|^2 \\ + \big((x_{4,1} - 8) + (x_{4,2} - 8)\big) + 
\left\|x_5 + \left(\begin{array}{c} 3 \\ 3 \end{array}\right)\right\|^4 \\ + \big((x_{6,1} - 10) + (x_{6,2} - 10)\big) + \big((x_{7,1} + 10) + (x_{7,2} + 10)\big) \\
+ \left\|x_8 + \left(\begin{array}{c} 7 \\ 0 \end{array}\right)\right\|^2 + \big((x_{9,1} - 6) + x_{9,2}\big) + \left\|x_{10} - \left(\begin{array}{c} 0 \\ 8 \end{array}\right)\right\|^4,
\end{multline}
where $x_{i,j}$ is the $j^{th}$ state of agent $i$ and the per-agent objectives can be discerned in the obvious way. The constraints on the agents are
\vspace{-0.25cm}
\begin{equation}
g(x) = \left(\begin{array}{c} \|x_1\|^2 + \|x_2\|^2 + \|x_3\|^2 - 10 \\
                              \|x_4\|^2 + \|x_5\|^2 + \|x_6\|^2 - 50 \\
                              \|x_7\|^2 + \|x_8\|^2 + \|x_9\|^2 - 50 \\  
                              x_{1,1}^2 + x_{5,1} + x_{10,1}^2 - 50 \\ 
                              x_{4,2}^2 + x_{7,1} + x_{9,2} - 20 \\ 
                              \|x_8\|^2 + \|x_6\|^2 - 30 \end{array}\right) \leq 0.
\end{equation}

Each agent was also constrained to lie in the box $X_i = [-10,10]\times[-10,10]$. The Lipschitz constants of $g$ were
computed to be $K^g_{1} = 39.82$ and $K^g_{2} = 56.71$. The Lipschitz constants for each $g_{x_i}$ are shown in Table \ref{tab:lipschitz1}.

\begin{table}[!htbp]
\centering
\begin{tabular}{|c|c|c|}
\hline
$i$ & $K^i_1$ & $K^i_2$ \\ \hline \hline 
$1$ & $4$ & $\sqrt{8}$ \\ \hline
$2$ & $2$ & $2$ \\ \hline
$3$ & $2$ & $2$ \\ \hline
$4$ & $2$ & $2$ \\ \hline 
$5$ & $2$ & $2$ \\ \hline 
$6$ & $4$ & $\sqrt{8}$ \\ \hline
$7$ & $2$ & $2$ \\ \hline 
$8$ & $4$ & $\sqrt{8}$ \\ \hline
$9$ & $2$ & $2$ \\ \hline
$10$& $2$ & $2$ \\ \hline
\end{tabular} 
\caption{Values of $K^i_1$ and $K^i_2$ for $g_{x_i}$, $i \in \{1, \ldots, 10\}$}
\label{tab:lipschitz1}
\end{table}

In both simulation runs below, the step-size rule discussed at the end of Section \ref{sec:conv} was used
with the values
\begin{equation}
\bar{\alpha} = 0.1, \,\, \bar{\gamma} = 0.01, \,\, c_1 = 0.3, \textnormal{ and } c_2 = 0.52,
\end{equation}
and all states and Kuhn-Tucker multipliers were initialized to zero, i.e., $x_i(0) = 0$ for all $i \in I$ and $\mu(0) = 0$.

\begin{figure}
\centering
\includegraphics[width=3.3in]{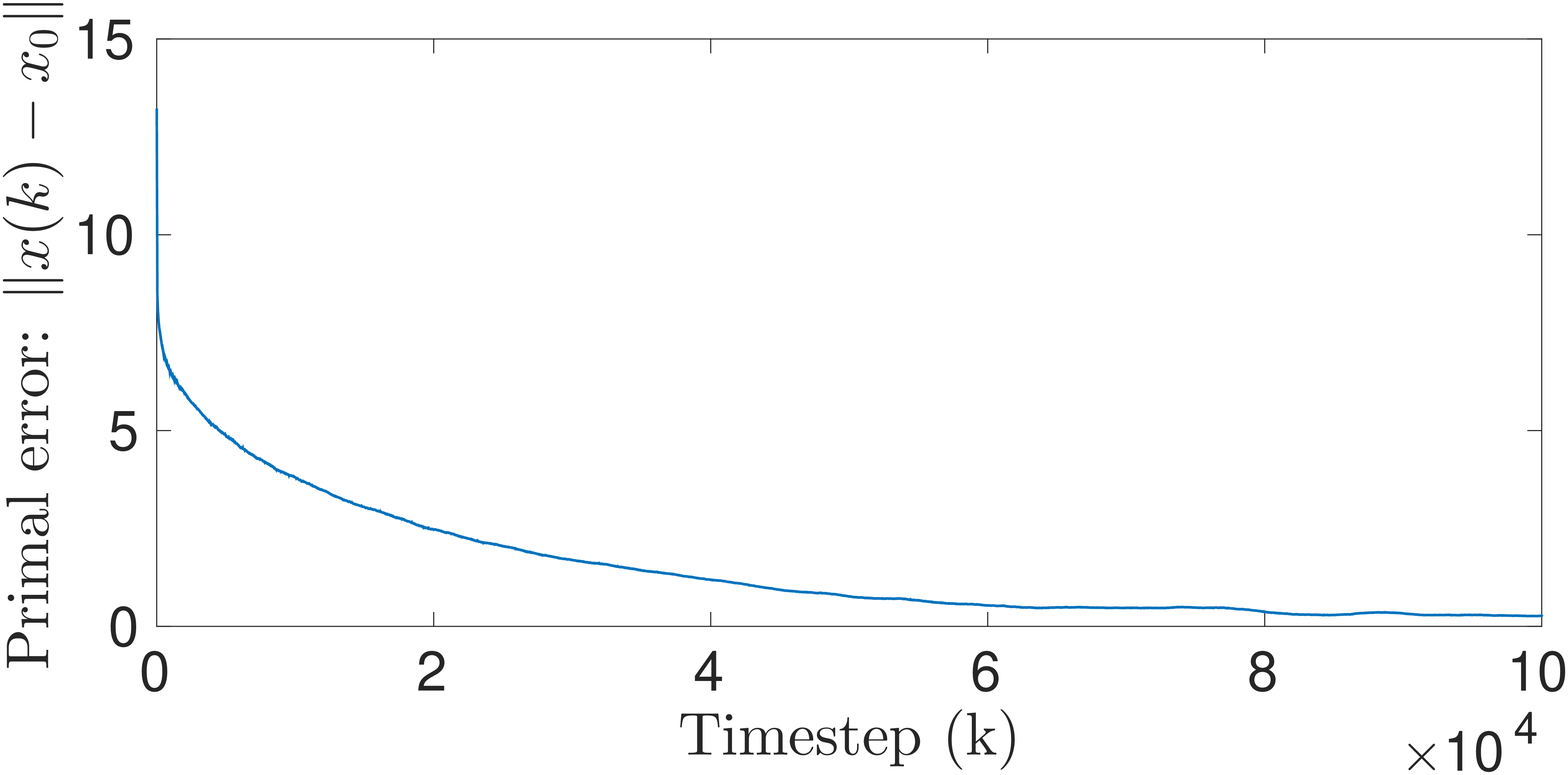}
\caption{\new{The values of $\|x(k) - x_0\|_2$ 
for $k = 1, \ldots, 100,000$ under $\epsilon$-differential privacy with Algorithm \ref{alg:stoch}.
The steady, monotone descent toward $x_0$ indicates numerical
convergence to $x_0$ in the presence of noise. }
}
\label{fig:epserrorx}
\end{figure}

\subsection{Simulation of $\epsilon$-differential privacy}
The adjacency parameter was chosen to be $B = 1$. 
The value $\epsilon = \ln 2$ was used for all systems. 
The distribution and variance of each entry of each noisy signal were 
computed and are listed in Table \ref{tab:lipschitz2} where we use the notation
for the PDF of a random scalar with the understanding that each entry
of the random matrices $w_i$ was generated using such a distribution.

\begin{table}[!htbp]
\centering
\begin{tabular}{|c|c|c|}
\hline
Noise & Distribution & Variance      \\ \hline \hline
$w_1$    & $\Lap(0, 5.771)$  & $66.60$ \\ \hline
$w_2$    & $\Lap(0, 2.885)$  & $16.65$ \\ \hline
$w_3$    & $\Lap(0, 2.885)$  & $16.65$ \\ \hline
$w_4$    & $\Lap(0, 2.885)$  & $16.65$ \\ \hline
$w_5$    & $\Lap(0, 2.885)$  & $16.65$ \\ \hline
$w_6$    & $\Lap(0, 5.771)$  & $66.60$ \\ \hline
$w_7$    & $\Lap(0, 2.885)$  & $16.65$ \\ \hline
$w_8$    & $\Lap(0, 5.771)$  & $66.60$ \\ \hline
$w_9$    & $\Lap(0, 2.885)$  & $16.65$ \\ \hline
$w_{10}$ & $\Lap(0, 2.885)$  & $16.65$ \\ \hline
\end{tabular} 
\caption{Noisy signals and their distributions for $\epsilon$-differential privacy}
\label{tab:lipschitz2}
\end{table}

\begin{figure}
\centering
\includegraphics[width=3.3in]{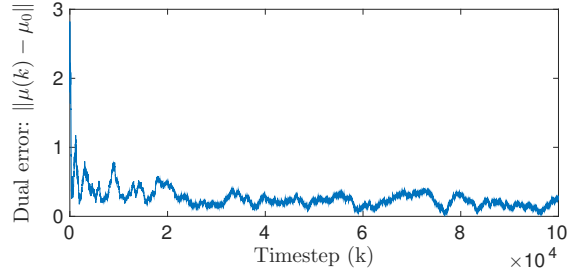}
\caption{\new{The values of $\|\mu(k) - \mu_0\|_2$ 
for $k = 1, \ldots, 100,000$ under $\epsilon$-differential privacy with Algorithm \ref{alg:stoch}.
Here we see an initial descent followed by a period of
oscillations as $\mu(k)$ approaches $\mu_0$. 
}
}
\label{fig:epserrory}
\end{figure}

The distribution for $w_g$ was $\Lap(0, 57.45)$ with variance $6.600 \cdot 10^3$.
Using this problem formulation, Algorithm \ref{alg:stoch}
was run for $100,000$ iterations. To show the behavior of the algorithm over time,
the least-norm saddle point of $L$, $z_0 = (x_0, \mu_0)$, was computed ahead of time and the values of $\|x(k) - x_0\|_2$ and 
$\|\mu(k) - \mu_0\|_2$ are shown in Figures \ref{fig:epserrorx} and \ref{fig:epserrory}, respectively, for $1 \leq k \leq 100,000$\footnote{Though the $1$-norm is used for
other aspects of $\epsilon$-differential privacy, we measure distance to $z_0$ using the $2$-norm to allow for meaningful visual
comparison of the plots corresponding to $\epsilon$-differential privacy in this subsection to those corresponding to $(\epsilon,\delta)$-differential privacy
in the next subsection.}.

In Figures \ref{fig:epserrorx} and \ref{fig:epserrory} we see a clear decreasing trend in both $\|x(k) - x_0\|_2$ and $\|\mu(k) - \mu_0\|_2$,
with the primal error appearing to be monotonically decreasing and the dual error oscillating while showing a general decreasing trend. 
The oscillations seen are expected given that
the variance of the noises added is constant while $G$ decreases in magnitude as
the saddle point $z_0$ is approached. In fact, it is known that descent will be achieved in a gradient method as long as the norm
of noise added to the gradient is less than the norm of the gradient itself \cite{zoutendijk70}. In light of this fact, the 
trends seen in Figures \ref{fig:epserrorx} and \ref{fig:epserrory} are not surprising because the gradients in $G$ in Algorithm \ref{alg:stoch} will have large
norms far from $z_0$, thereby allowing them to ``overpower'' the noise added, while close to $z_0$ their
norms will be smaller and the noise can dominate, causing increases in the distance to $z_0$ at some timesteps. Of course,
$z(k) \to z_0$ asymptotically because these increases in $\|z(k) - z_0\|_2$ average out over very long periods of time.

The initial error values here were
\begin{equation}
\|x(0) - x_0\|_2 = 13.19 \textnormal{ and } \|\mu(0) - \mu_0\|_2 = 2.169,
\end{equation}
\new{
And in this run the final error values were
\begin{equation}
\|x(100,000) - x_0\|_2 = 0.2706 \textnormal{ and } \|\mu(100,000) - \mu_0\|_2 = 0.2842,
\end{equation}
with these values after half of the total runtime being
\begin{equation}
\|x(50,000) - x_0\|_2 = 0.7658 \textnormal{ and } \|\mu(50,000) - \mu_0\|_2 = 0.2225.
\end{equation}
}

\vspace{-0.3cm}

These values confirm what can be seen visually in Figures \ref{fig:epserrorx} and \ref{fig:epserrory}: 
shorter runtimes than $100,000$ timesteps can be used while ending at a reasonable distance from $z_0$ and, in light of the
large variances of some noises present, reasonable numbers of iterations produce an approach toward $z_0$ that would be 
useful in many applications. 

\begin{figure}
\centering
\includegraphics[width=3.3in]{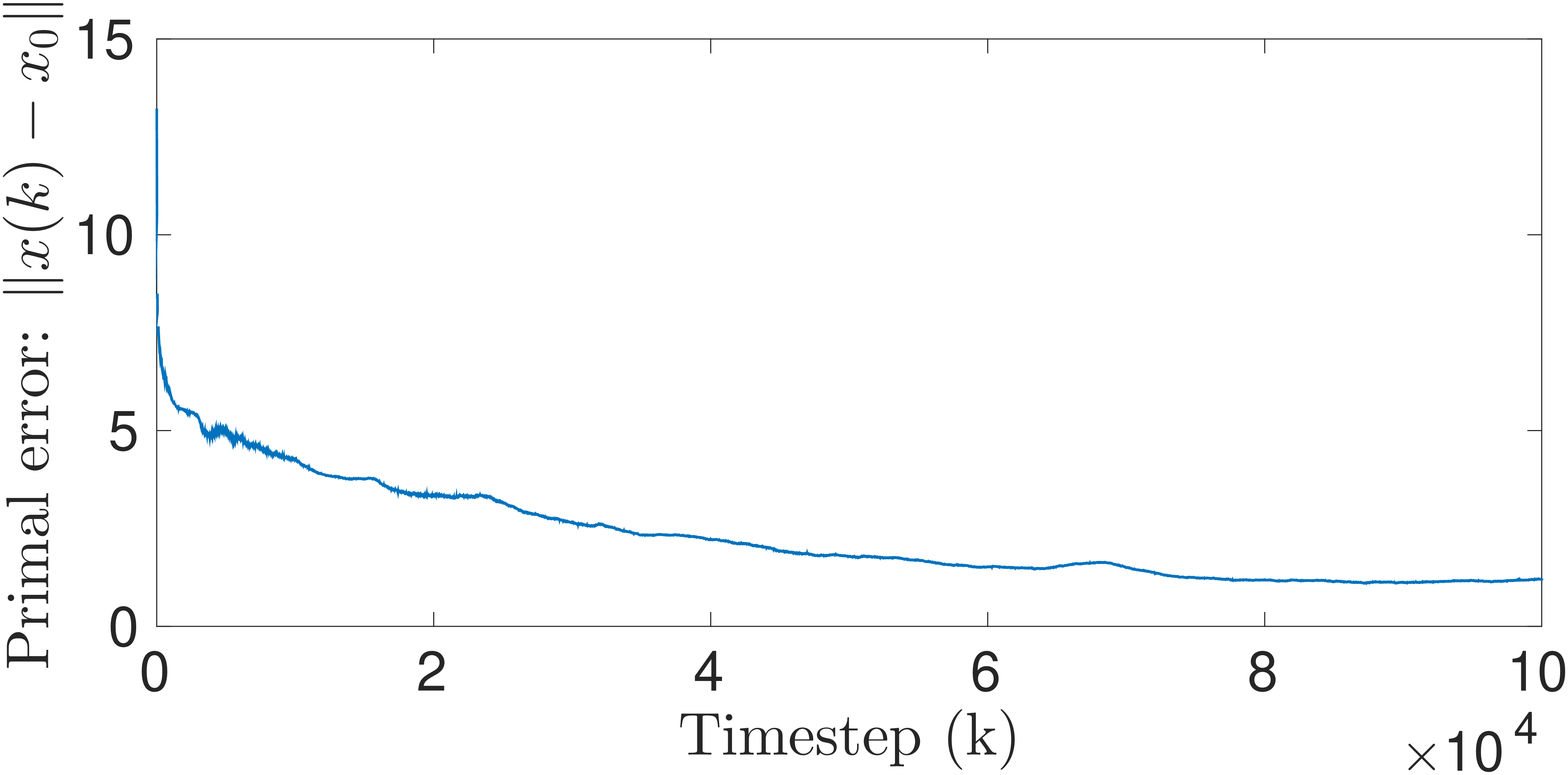}
\caption{\new{The values of $\|x(k) - x_0\|_2$ for $k = 1, \ldots, 100,000$ under $(\epsilon, \delta)$-differential privacy with Algorithm \ref{alg:stoch}.
The rapid descent toward $x_0$ and clear decreasing trend thereafter indicate numerical convergence to $x_0$ in the presence of noise. }
}
\label{fig:epsderrorx}
\end{figure}

\subsection{Simulation of $(\epsilon, \delta)$-differential privacy}
In this case the adjacency parameter was chosen to be $B > 1$. 
The values $\epsilon = \ln 2$ and $\delta = 0.01$ were used for all systems,
giving $\kappa(\delta, \epsilon) = 3.559$. Using 
this privacy policy, the distribution and variance of each noisy signal were computed and
are listed in Table \ref{tab:lipschitz3}. 

\begin{table}[!htbp]
\centering
\begin{tabular}{|c|c|c|}
\hline
Noise & Distribution & Variance      \\ \hline \hline
$w_1$    & $\mathcal{N}(0,101.3)$ & $101.3$ \\[2pt] \hline
$w_2$    & $\mathcal{N}(0,50.66)$ & $50.66$ \\[2pt] \hline
$w_3$    & $\mathcal{N}(0,50.66)$ & $50.66$ \\[2pt] \hline
$w_4$    & $\mathcal{N}(0,50.66)$ & $50.66$ \\[2pt] \hline
$w_5$    & $\mathcal{N}(0,50.66)$ & $50.66$ \\[2pt] \hline
$w_6$    & $\mathcal{N}(0,101.3)$ & $101.3$ \\[2pt] \hline
$w_7$    & $\mathcal{N}(0,50.66)$ & $50.66$ \\[2pt] \hline
$w_8$    & $\mathcal{N}(0,101.3)$ & $101.3$ \\[2pt] \hline
$w_9$    & $\mathcal{N}(0,50.66)$ & $50.66$ \\[2pt] \hline
$w_{10}$ & $\mathcal{N}(0,50.66)$ & $50.66$ \\ \hline
\end{tabular} 
\caption{Noisy signals and their distributions for $(\epsilon,\delta)$-differential privacy}
\label{tab:lipschitz3}
\end{table}


The distribution for $w_g$ was $\mathcal{N}(0_{6 \times 1}, 4.073 \cdot 10^4 I_{6 \times 6})$ with variance $4.073 \cdot 10^4$. 
\new{
In Table \ref{tab:lipschitz3} we record the distribution of each entry of the matrices $w_i$, $i \in I$, with the understanding
that each $w_i$ has i.i.d. entries. 
}

\begin{figure}
\centering
\includegraphics[width=3.3in]{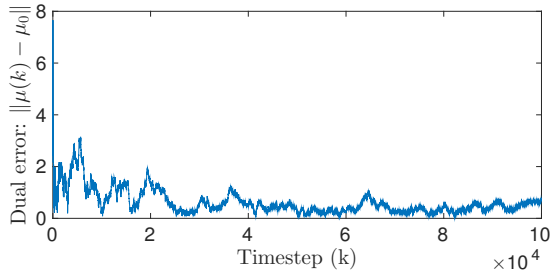}
\caption{\new{The values of $\|\mu(k) - \mu_0\|_2$ for $k = 1, \ldots, 100,000$ under $(\epsilon, \delta)$-differential privacy with Algorithm \ref{alg:stoch}.
The initial approach toward $\mu_0$ and oscillations in distance beyond that point indicate numerical convergence to $\mu_0$ 
when noise is added for differential privacy.}
}
\label{fig:epsderrory}
\end{figure}

Using this problem formulation, Algorithm \ref{alg:stoch}
was run for $100,000$ iterations and the values of $\|x(k) - x_0\|_2$ and
$\|\mu(k) - \mu_0\|_2$ for $1 \leq k \leq 100,000$ are plotted in Figures 
\ref{fig:epsderrorx} and \ref{fig:epsderrory}, respectively. 
In Figures \ref{fig:epsderrorx} and \ref{fig:epsderrory} we see a similar trend to Figures \ref{fig:epserrorx} and \ref{fig:epserrory}: 
steady, monotone decrease in the primal error, and general decreases in the dual error with noticeable oscillations present. 

The initial error values for this run were
\begin{equation}
\|x(0) - x_0\|_2 = 13.19 \textnormal{ and } \|\mu(0) - \mu_0\|_2 = 2.169.
\end{equation}
\new{
The final error values here were
\begin{equation}
\|x(100,000) - x_0\|_2 = 1.1965 \textnormal{ and } \|\mu(100,000) - \mu_0\|_2 = 0.7413,
\end{equation}
while after half of the total timesteps taken these values were
\begin{equation}
\|x(50,000) - x_0\|_2 = 1.7857 \textnormal{ and } \|\mu(50,000) - \mu_0\|_2 = 0.2500,
\end{equation}
indicating a rapid initial descent towards $z_0$ and close proximity to it thereafter. 
}

Both simulation examples show a rapid decrease in the distance from $z(k)$ to $z_0$. Such a rapid decrease lends itself to use of this algorithm
in practical applications because it allows for useful improvements to be made in the value of $f$ in a reasonable runtime while respecting the
set and functional constraints of the problem. 
The theoretical and simulation results presented here demonstrate the utility of the iterative Tikhonov regularization, even in the presence of noise with large variance.
This robustness is further supported by the simulation results in \cite{hale15} and demonstrates that in a practical setting strong, quantifiable guarantees of privacy can be achieved 
while providing useful convergence guarantees in the optimization problem of interest. 
Critical to the success of these numerical results is all noise being zero mean,
 and it is a feature of differential privacy that zero mean noise is effective at protecting sensitive information. 

\section{Conclusion}
A differentially private optimization algorithm for teams of many agents coordinated by a central cloud computer
was presented. This problem was treated as a stochastic 
variational inequality and solved using a Tikhonov-regularized
Goldstein-Levitin-Polyak iteration. Its convergence was shown for both
$\epsilon$- and $(\epsilon, \delta)$-differential privacy and numerical
convergence of the algorithm was shown in simulation, demonstrating
the ability to arrive at a collective decision while maintaining
privacy for the users involved in making it.

\bibliographystyle{IEEEtran}{}
\bibliography{sources}

\end{document}